\documentclass[12pt, a4paper]{amsart}

\usepackage{latexsym, amssymb,amsmath, amsthm,amsfonts}
\usepackage{hyperref,url}

\newcommand{\FF}{\mathbb{F}}

\newcommand{\ZZ}{\mathbb{Z}}
\newcommand{\kk}{\Bbbk}
\newcommand{\kv}{{\kk[V]}}

\newcommand{\kvgbar}{\overline{\kk}[\overline{V}]^G}
\newcommand{\vbar}{\overline{V}}
\newcommand{\kbar}{\overline{\kk}}
\newcommand{\kvbar}{\overline{\kk}[\overline{V}]}
\newcommand{\kvg}{{\kk[V]^G}}
\newcommand{\kvgplus}{{\kk[V]^G_+}} 

\newcommand{\done}{\hfill $\triangleleft$}

\def\SL2{\operatorname{SL}_{2}(K)}

\def\GL2{\operatorname{GL}_{2}(K)}

\def\INVSL2{$K[V]^{operatorname{SL}_{2}(K)}$}
\def\INVSO2{$K[V]^{operatorname{SO}_{2}(K)}$}
\def\INVGL2{$K[V]^{operatorname{GL}_{2}(K)}$}

\def\depth{\operatorname{depth}}
\def\height{\operatorname{ht}}
\def\cmdef{\operatorname{cmdef}}

\def\Hom{\operatorname{Hom}}
\def\GL{\operatorname{GL}}
\def\SL{\operatorname{SL}}

\def\Z{\mathbb{Z}}

\def\Ann{\operatorname{Ann}}

\def\codim{\operatorname{codim}}

\newtheorem{Lemma}{Lemma}[section]
\newtheorem{Theorem}[Lemma]{Theorem}
\newtheorem{Corollary}[Lemma]{Corollary}

\newtheorem{Prop}[Lemma]{Proposition}

\theoremstyle{definition}

\theoremstyle{remark}
  \newtheorem{rem}[Lemma]{Remark}
\newtheorem{eg}[Lemma]{Example}

\newtheoremstyle{Acknowledgements}% name
  {}% {\topsep}%      Space above
    {}% {\topsep}%      Space below
     {}%         Body font
     {}%         Indent amount (empty = no indent, \parindent = para indent)
    {\bfseries}% Thm head font
    {}%        Punctuation after thm head
     {.5em}%     Space after thm head: " " = normal interword space;
     {\thmname{#1}\thmnumber{ }\thmnote{ (#3)}}%         Thm head spec (can be left empty, meaning `normal')
\theoremstyle{Acknowledgements}
\newtheorem{ack}{Acknowledgements.}

\title[ The C-M property of separating invariants]
{The Cohen-Macaulay property  of \\ separating invariants of finite groups}

\author{ Emilie Dufresne}
\address{Mathematics Center Heidelberg (MATCH) \\
Ruprecht-Karls Universit\"at Heidelberg\\
Im Neuenheimer Feld 368\\
69120 Heidelberg, Germany}
\email{emilie.dufresne@iwr.uni-heidelberg.de}

\author{Jonathan Elmer}
\address{Lehrstuhl A f\"ur Mathematik\\
RWTH Aachen\\
52062 Aachen, Germany}
\email{jonathan.elmer@mathA.rwth-aachen.de}

\author{ Martin Kohls}
\address{Technische Universit\"at M\"unchen \\
 Zentrum Mathematik-M11\\
Boltzmannstrasse 3\\
 85748 Garching, Germany}
\email{kohls@ma.tum.de}

\date{\today}

\begin{document}

%%%%%%%%%%%%%%%%%%%%%%%%%%%%%%%%%%%%%%%%%%%%%%%%%%%%%%%%%%%%%%

\begin{abstract}{In the case of finite groups, a separating algebra is a subalgebra of the ring
of invariants which separates the orbits. Although separating algebras are
often better behaved than the ring of invariants, we show that many of the criteria which imply the ring of invariants is non Cohen-Macaulay
actually imply that no graded separating algebra is Cohen-Macaulay. For
example, we show that, over a field of positive characteristic $p$, given
sufficiently many copies of a faithful modular representation, no graded separating algebra is Cohen-Macaulay. Furthermore, we show that, for a $p$-group, the existence of a Cohen-Macaulay graded separating algebra implies the group is generated by bireflections. Additionally, we give an example which shows that 
Cohen-Macaulay separating algebras can occur when the ring of invariants is not Cohen-Macaulay.
}\end{abstract}

%%%%%%%%%%%%%%%%%%%%%%%%%%%%%%%%%%%%%%%%%%%%%%%%%%%%%%%%%%%%%%

\maketitle

%%%%%%%%%%%%%%%%%%%%%%%%%%%%%%%%%%%%%%%%%%%%%%%%%%%%%%%
%%%%%%%%%%%%%%%%%%%%%%%%%%%%%%%%%%%%%%%%%%%%%%%%%%%%%%%%%
%%%     SECTION:       Introduction                        %%%%%%%%%%%%%%%%%%%%%%%%%%%%%%%%%%%%%%%%%%%%%%%%%%%%%%%%%%%%%%%%%%%%%%%%%%%%%%%%%%%%%%%%%%%%%%%%%%%%%%%%%%%%%%%%%%%%%%%%%%%%%%%%%%%%%%%%%%%%%%%%%%%%%%%%%%%%%%%%%%%%%%%%%%%%%%%%%%%%%%%%%%%%%%%%%%%%%%%%%%%%%%%%%%%%%%%%%%%%%%%%%%%%%%%%%%%%%%%%%%%%%%%%%%%%%%%%%%%%%%%%%%%%%%%%%%%%%%

\section{Introduction}

Let $G$ be  a finite group, and let $V$ be a finite dimensional representation of $G$  over a field $\kk$ of characteristic $p\geq 0$. We say $V$ is a modular representation if $p$ divides $|G|$. We write $\kv$ for the symmetric algebra $S(V^*)$ on the vector space dual of $V$. It is a polynomial ring with the standard grading.  The action of $G$ on $V$ induces an action on $\kv$: for  $f\in V^{*}$ and $v\in V$,  the action of  $\sigma\in G$  is given by $(\sigma\cdot f)(v)=f(\sigma^{-1}\cdot v)$. The ring of invariants, denoted $\kvg$, is the ring formed by the elements of $\kv$ fixed by $G$. Since the group action preserves degree, $\kvg$ is a graded subalgebra of $\kv$.

Let $\kbar$ be an algebraic closure of the field $\kk$, and let $\vbar=\kbar\otimes_{\kk}V$. As $\kv\subseteq\kbar\otimes_{\kk}\kv\cong\kvbar $,  any  $f$ in $\kv$ can be considered as a function $\vbar\rightarrow{\kbar}$. The action of $G$ on $V$ extends to an action of $G$ on $\vbar$, and so $\kvg\subseteq\kvgbar$. 

By definition, elements of $\kvg$ are constant on $G$-orbits. Accordingly, if an invariant $f$ takes distinct values on elements $u,v\in\vbar$, then these elements belong to distinct orbits, and we say $f$
\emph{separates} $u$ and~$v$. A \emph{geometric separating set} is then a set of invariants which separates exactly the same points of $\vbar$ as the whole ring of invariants. As $G$ is finite, the ring of invariants separates orbits in $\vbar$ \cite[Lemma~2.1]{DraismaSeparating}. Hence, a geometric separating set is a set of invariants which separates the orbits of the $G$-action on $\vbar$.

The study of separating invariants, which has become quite popular in recent
years, was initiated by Derksen and Kemper (\cite[Section 2.3]{DerksenKemper}
and \cite{KemperCompRed}). They defined a \emph{separating set} as a set separating the same orbits in $V$ as the whole ring of invariants \cite[Definition~2.3.8]{DerksenKemper}.
For $\kk$ algebraically closed, geometric separating sets coincide with separating sets, but in general, a separating set is not always a geometric separating set (see Example~\ref{C4perm}). On the other hand, $A\subseteq \kvg$ is a geometric separating algebra if and only if  $\kbar\otimes_{\kk}A\subseteq \kvgbar$ is a separating algebra. Moreover, if $A$ is graded, so is $\kbar\otimes_\kk A$, and the two rings have the same dimension and depth. Furthermore, the extension $A\subseteq\kvg$ is integral if and only if $\kbar\otimes_\kk A\subseteq\kvgbar$ is integral. Thus, it often suffices to write proofs for separating sets over algebraically closed fields. Note that this works only because we are interested in geometric separating algebras.

Many defects of invariant rings disappear when one considers separating
invariants. The ring of invariants is not always finitely generated (for
non-reductive groups) \cite[Example~2.1.4, due to Nagata]{DerksenKemper}, but
there always exist finite geometric separating sets
\cite[Theorem~2.3.15]{DerksenKemper}. Over algebraically closed fields, there
is an upper bound on the size of minimal separating sets, depending only on the dimension of the representation \cite[Proposition~5.1.1]{DufresnePhd}. Polarization, a classical method for obtaining vector invariants in characteristic zero, extends, for separating invariants, to all characteristics \cite{DraismaSeparating, DomokosSep}. For $G$ finite, the Noether bound holds for separating invariants in all characteristics: although they may not generate the ring of invariants, the invariants of degree at most $|G|$ always form a geometric separating set \cite[Section~3.9]{DerksenKemper}.

Graded separating algebras are very closely related to the ring of invariants:
\begin{Prop}\label{FiniteGenerationOfGeometric}
If $A\subseteq \kvg$ is a graded geometric separating algebra, then  $A\subseteq \kvg$ is an integral extension, and $A$ is a finitely generated $\kk$-algebra. 
\begin{proof}
We may assume $\kk$ is algebraically closed. By
\cite[Lemma~1.3]{KemperCompRed}, the extension $A\subseteq\kvg$ is
integral. The finite generation of $A$ as a $\kk$-algebra then follows from that
of $\kvg$ by Newstead \cite[p. 52, (II)]{Newstead}, or in the manner of \cite[Proof of Theorem 1.3.1]{Benson}.
\end{proof}\end{Prop}

\begin{Prop}\label{pRootClosure}
Suppose $p> 0$. If $A\subseteq\kvg$ is a graded subalgebra, then $A$ is a geometric separating algebra if and only if $\kvg$ is the purely inseparable closure of $A$ in $\kv$, that is,
\[
\kvg=\{f\in \kv \mid \quad \textrm{for some}~m,~f^{p^{m}}\in A\}.\] 
\begin{proof}
For $\kk$ algebraically closed, see \cite[Remark 1.3]{DerksenKemperAlgebraic}. The proof is an application of a result of van der Kallen \cite[Sublemma
A.5.1]{VanDerKallen} (see the extended proof in \cite{VanDerKallenOnline}). A
variation of the same argument first appeared in \cite[Theorem 6]{GrosshansVectorInvariants}.  For $f\in \kvg$ and $\kk$ arbitrary, we have
that $f^{p^{m}}\in \kbar\otimes_{\kk} A$ for some $m$. The stability of the
rank of matrices under field extensions implies $f^{p^{m}}\in A$.
\end{proof}
\end{Prop}

\begin{rem}
 Propositions~\ref{FiniteGenerationOfGeometric} and \ref{pRootClosure} hold for rational representations of reductive groups.
\end{rem}

Kemper~\cite{KemperCompRed} exploited this close relationship to compute the
invariants of reductive groups in positive characteristic. On the other hand,
Dufresne~\cite{DufresneSeparating} showed that the existence of
polynomial or complete intersection separating algebras imposes strong
conditions on the representation. The present paper is in the latter vein. We show that, in many instances, conditions which ensure that the ring of invariants is non Cohen-Macaulay, in fact imply that no graded geometric separating algebra is Cohen-Macaulay. We thus provide a (partial) negative answer to Kemper who asked if Cohen-Macaulay separating algebras should always exist~\cite{KemperSeparating}. Notably, we show:

\begin{Theorem}\label{kcopiesTheorem}
If $V$ is faithful and modular, then there exists $r\!\geq 1$ such that, for all
$k$, every graded geometric separating algebra in $\kk[V^{\oplus k}]^{G}$ has
Cohen-Macaulay defect at least $k-r-1$. In particular, for $k>r+1$, no graded geometric separating algebra in $\kk[V^{\oplus k}]^{G}$ is Cohen-Macaulay.
\end{Theorem}

An element $\sigma$ of $G$ acts as a
\emph{bireflection} on $V$ if its fixed space is of codimension at most $2$ in $V$. 

\begin{Theorem}\label{thm1}
Let $G$ be a $p$-group. If there exists a graded geometric separating algebra
in $\kvg$ which is Cohen-Macaulay, then $G$ is generated by elements acting as bireflections.
\end{Theorem}

Theorem~\ref{thm1} fits well with \cite[Theorem~1.3]{DufresneSeparating}. In
the important special case of $p$-groups, we obtain that $G$ is generated by
bireflections from a much weaker
hypothesis: the existence of a Cohen-Macaulay rather than a complete
intersection graded geometric separating algebra. This mirrors the situation for invariant rings
(\cite[Corollary~3.7]{KemperOnCM} and \cite[Theorem~A]{vk-kw:flgrici}). Example~\ref{refnEg} shows that not only the converse of Theorem~\ref{thm1}, but also the converses of \cite[Theorem~1.1 and Theorem~1.3]{DufresneSeparating}, are not true.

In Section~\ref{FirstCohom}, we extend the methods introduced in \cite{KemperOnCM} to prove our main results. Section~\ref{A4} concentrates on  the alternating group $A_4$. We conclude in Section \ref{remarks} with a discussion of the general situation and examples which show that the depth of graded geometric separating algebras can be both larger and smaller than that of the corresponding invariant ring.

\begin{ack}
This paper was prepared during visits of the first and second authors to TU
M\"unchen, and of the second and third authors to Universit\"at Heidelberg. We
 thank Gregor Kemper, B. Heinrich Matzat, and MATCH for making these
visits possible. Finally, we thank the anonymous referees for helpful
suggestions, in particular for pointing out an error in our original proof of Lemma \ref{trivcoeff}.
\end{ack}

%%%%%%%%%%%%%%%%%%%%%%%%%%%%%%%%%%%%%%%%%%%%%%%%%%%%%%%%%%%%%%%%%%%%%%%%%%%%%%%%%%%%%%%%%%%%%%%%%%%%%%%%%%%%%%%%%%%%%%%%%%%%%%%%%%%%%%%%%%%%%%%%%%%%%%%%%%%%%%%%%%%%%%%%%%
%%%%    Section:  The Cohen-Macaulay Defect of Separating Algebras       %%%%%%%%%%%%%%%%%%%%%%%%%%%%%%%%%%%%%%%%%%%%%%%%%%%%%%%%%%%%%%%%%%%%%%%%%%%%%%%%%%%%%%%%%%%%%%%%%%%%%%%%%%%%%%%%%%%%%%%%%%%%%%%%%%%%%%%%%%%%%%%%%%%%%%%%%%%%%%%%%%%%%%%%%%%%%%%%%%%%%%%%%%%%%%%

\section{ The Cohen-Macaulay Defect of Separating Algebras}\label{FirstCohom}

Let $A\subseteq \kv$ be a finitely generated graded subalgebra, and let
$A_{+}$ denote its maximal homogeneous ideal. Homogeneous elements
$a_{1},\ldots,a_{k}$ in $A_{+}$ form a \emph{partial homogeneous system of
  parameters (phsop)} if they generate an ideal of height $k$ in $A$. If
additionally $k=\dim A$, then they form a \emph{homogeneous system of
  parameters (hsop)}. Noether's normalization theorem guarantees that a hsop
always exists. If, for $i=1,\ldots,k$, the element $a_{i}$ is not a zero
divisor in $A/(a_{1},\ldots,a_{i-1})A$, then the elements $a_{1},\ldots,a_{k}$
form a \emph{regular sequence}. Every regular sequence is a phsop. We say $A$
is \emph{Cohen-Macaulay} when every phsop is a regular sequence. The depth of a
homogeneous ideal $I\subseteq A_+$, written $\depth_A(I)$, is the maximal length of a
regular sequence in $I$. Note that the height of $I$, $\height(I)$, is equal to the maximal length of a phsop in $I$. We write $\depth (A):=\depth_A(A_{+})$, and define the \emph{Cohen-Macaulay defect} of $A$ to be $\cmdef (A):=\dim A-\depth (A)$. Thus, $A$ is Cohen-Macaulay precisely when $\cmdef A=0$.

In Theorem \ref{MainTheorem}, we relate the Cohen-Macaulay defect of graded geometric separating algebras  to the $n$-th cohomology group $H^{n}(G,\kvbar)$. For the theory of these
 groups for arbitrary $n$, we refer to \cite{BensonCohomology1,Weibel}. Since we have $n=1$ in most applications of this theorem, we construct this group now.

 Let $W$ be a representation of the group $G$ over the field $\kk$. A
 \emph{1-cocycle} is a map $g:G\rightarrow W,\,\sigma\mapsto g_{\sigma}$
 such that $g_{\sigma\tau}=\sigma g_{\tau}+g_{\sigma}$, for all
 $\sigma,\tau\in G$. We write $Z^{1}(G,W)$ for the additive group of all
 1-cocycles. For each $w\in W$, the map given by $\sigma\mapsto
 (\sigma-1)w:=\sigma w- w$ is a 1-cocycle, which is called a \emph{1-coboundary}. The coboundaries form the subgroup
 $B^{1}(G,W)$ of $Z^{1}(G,W)$. The \emph{first cohomology group of $G$ with
 coefficients in $W$} is the quotient group
 $H^{1}(G,W):=Z^{1}(G,W)/B^{1}(G,W)$. A cocycle $g$ is \emph{nontrivial} if
 and only if its \emph{cohomology class} $g+B^{1}(G,W)$ is nonzero in
 $H^{1}(G,W)$. Similar definitions hold for $H^{n}(G,W)$. We will sometimes abuse
notation by using the same symbols for cocycles and cohomology classes. 

The group $H^{n}(G,\kv)$ has a natural graded $\kvg$-module structure. For
a homogeneous $g\in H^{n}(G,\kv)$,  its  annihilator in
$\kvg$, 
\[\Ann_{\kvg}(g):=\{a\in \kvg: ag=0\in H^{n}(G,\kv)\},\]
is a homogeneous ideal. If $p>0$, the $m$-fold Frobenius homomorphism
$\kv\rightarrow\kv,\, f\mapsto f^{p^{m}}$, induces a map $H^{n}(G,\kv)\rightarrow H^{n}(G,\kv)$. This map is a homomorphism of abelian
groups, but not of $\kk$-vector spaces. We write
$g^{p^{m}}$ for the image of an element $g\in H^{n}(G,\kv)$ under this map.
In particular, for $g\in H^{1}(G,\kv)$, the cohomology class $g^{p^{m}}\in H^{1}(G,\kv)$ is given by the cocycle $\sigma\mapsto {(g_\sigma)}^{p^m}$\!\!\!.

Over fields of characteristic zero, and in the non-modular case in general,
the ring of invariants is always Cohen-Macaulay \cite{HochsterEagon}. In
particular, there will always be a Cohen-Macaulay geometric separating
algebra. Accordingly, from now on we assume $V$ is modular.

Our most general statement generalizes \cite[Corollary
1.6]{KemperOnCM} (see also \cite[Proposition 6]{KemperLinRed}) to the case of
separating algebras:

%%%%%% Main theorem %----------------------

\begin{Theorem}\label{MainTheorem}
Let $n\geq 1$ be the smallest integer such that there exists a homogeneous $g\in
H^{n}(G,\kvbar)$ such that $g^{p^{m}}$ is nonzero for every $m\ge 0$.
If $A$ is a graded geometric separating algebra in $\kvg$, then $I:=\Ann_{A}(g)$, the
annihilation ideal of $g$ in $A$, has depth at most $n+1$. Furthermore, if
$\Ann_{\kvgbar}(g)$ has height $k$, then $A$ has Cohen-Macaulay defect at least $k-n-1$.
\end{Theorem}

Since the case $n=1$ suffices for most of our applications, we give an
additional more elementary proof of the first part. Without loss of generality, we assume $\kk=\kbar$ in both arguments.

\begin{proof}[Proof of the case $n=1$.] 
 By Proposition~\ref{pRootClosure}, there exists a $p$-power $q$ such that
  $(\kvg)^{q}\subseteq A$. Suppose, for a contradiction, that $\depth_{A}(I)$
  is at least $3$. Hence, there exists an $A$-regular
sequence $a_{1},a_{2},a_{3}$ in $I$. Since $a_{i}g=0$,
there are $b_{1},b_{2},b_{3}\in \kv$ such that
\[
a_{i}g_{\sigma}=(\sigma-1)b_{i},\quad \text{for all }\sigma\in G,\,\,i=1,2,3.
\]
Set $u_{ij}:=a_{i}b_{j}-a_{j}b_{i}$, for $1\le i< j \le 3$. 
For all $i,j$, $u_{ij}$ is invariant, and so
$u_{ij}^{q}$ belongs to $A$. Since  $a_{1}^{q},a_{2}^{q},a_{3}^{q}$ forms an $A$-regular sequence  \cite[Corollary 17.8 (a)]{Eisenbud}, and since \renewcommand{\arraystretch}{0.7}
\renewcommand{\arraycolsep}{2pt}
\[
a_{1}^{q}u_{23}^{q}-a_{2}^{q}u_{13}^{q}+a_{3}^{q}u_{12}^{q}=\left|\begin{array}{ccc}
    a_{1}&a_{2}&a_{3}\\ a_{1}&a_{2}&a_{3}\\
    b_{1}&b_{2}&b_{3}\end{array}\right|^{q}=0,
\]
it follows that $u_{12}^{q}\in(a_{1}^{q},a_{2}^{q})A$. Thus, there exist
$f_{1},f_{2} \in A$ such that
\[
u_{12}^{q}=a_{1}^{q}b_{2}^{q}-a_{2}^{q}b_{1}^{q}=f_{1}a_{1}^{q}+f_{2}a_{2}^{q}.
\]
As $a_{1},a_{2}$ is a phsop in $A$, it is also a phsop in its integral
extension $\kv$, and thus $a_{1},a_{2}$ are coprime in $\kv$. From
$a_{1}^{q}(b_{2}^{q}-f_{1})=a_{2}^{q}(f_{2}+b_{1}^{q})$, it follows that
$a_{1}^{q}$ divides $f_{2}+b_{1}^{q}$ in $\kv$. Therefore, there exists $h\in \kv$
such that $a_{1}^{q}h=f_{2}+b_{1}^{q}$. Hence, for every $\sigma\in G$, we have
\[
a_{1}^{q}(\sigma-1)h=((\sigma-1)b_{1})^{q}=(a_{1}g_{\sigma})^{q}=a_{1}^{q}g_{\sigma}^{q},
\]
that is, $g_{\sigma}^{q}=(\sigma-1)h$, a contradiction since $g^{q}$ is
nonzero. Thus, $\depth_{A}(I)\leq 2$.
\end{proof}

%--------------------------------------------------------------
% proof of general case
%-------------------------------------------------------------

\begin{proof}[Proof of the general case.]
For some $m\geq 0$, we have $(\kk [V]^G)^{p^m} \subseteq A$. For each $i$,
$H^i(G,\kk [V])$ is finitely generated as a $\kk [V]^G$-module. Therefore,
for each $0<i<n$, there exists some $m(i)$ such that $\alpha^{p^{m(i)}} = 0$ for
all $\alpha \in H^i(G,\kk [V])$. Set $q:= p^{m'}$, where $m'$ is the maximum of $m, m(1), m(2), \ldots, m(n-1)$. Assume, for a contradiction, that $\depth_A(I) \geq n+2$, and let $a_1, a_2,
\ldots , a_{n+2}$ be an $A$-regular sequence in~$I$. Consider the bar
resolution of $\Z$ as a ${\mathbb \Z}G$-module:
$$\ldots \stackrel{\partial_{n+1}}{\rightarrow} X_n \stackrel{\partial_n}{\rightarrow} X_{n-1} \stackrel{\partial_{n-1}}{\rightarrow} \ldots \stackrel{\partial_1}{\rightarrow} X_0 \stackrel{\partial_0}{\rightarrow} \mathbb{Z} \rightarrow 0,$$
where $X_0 = \mathbb{Z}G$,  and $\partial_0$ is the augmentation
map. The cohomology class $g$ is represented by a  $u_{\{1\}}\in
\Hom_{{\mathbb Z}G}(X_{n},\kv)$ such that $u_{\{1\}} \partial_{n+1} = 0$, that
is, $u_{\{1\}}\in Z^{n}(G,\kv)$
(the notation will become clear later). For each $i$, $a_i g = 0$, so there
is $h_{\{i\}} \in \Hom_{{\mathbb Z}G}(X_{n-1},\kv)$ such that $a_i
u_{\{1\}} = h_{\{i\}}\partial_n$. We next define, for each $2\le r \le n+1$, and for each ordered
$r$-subset $J:= \{j(1), j(2), \ldots ,j(r)\}$ of $\{1, \ldots, n+2\}$, a homomorphism
$u_J \in \Hom_{\mathbb{Z}G}(X_{n-r+1},\kk [V])$ such that $u_J \partial_{n-r+2} =
0$, that is, $u_{J}\in Z^{n-r+1}(G,\kv)$. For $2\le r\le n$, the definition of $q$ implies $u_{J}^{q}\in B^{n-r+1}(G,\kv)$,
so there is a map $h_J \in \Hom_{\mathbb{Z}G}(X_{n-r}, \kk [V])$ satisfying
 $u_J^q = h_J \partial_{n-r+1}$. We now define $u_{J}$, and thus $h_{J}$, by
 induction on $r$:
\begin{equation}\label{DefNOfuJ}
u_J := \sum_{i=1}^r (-1)^{i+1}a^{q^{r-2}}_{j(i)} h_{J \setminus \{j(i)\}}.
\end{equation}
 Next, we show that $u_J \partial_{n-r+2}=0$ for $2\le r \le n+1$. For $r=2$, we get
\begin{eqnarray*}
u_{\{j(1),j(2)\}}\partial_{n}&=&a_{j(1)}h_{\{j(2)\}}\partial_{n}-a_{j(2)}h_{\{j(1)\}}\partial_{n}\\&=&a_{j(1)}a_{j(2)}u_{\{1\}}-a_{j(2)}a_{j(1)}u_{\{1\}}=0.
\end{eqnarray*}
For $2<r\le n+1$, we obtain similarly:
\begin{equation}\label{uJdnr2}
u_J \partial_{n-r+2} = \sum_{i=1}^{r}(-1)^{i+1} a^{q^{r-2}}_{j(i)} u^q_{J
  \setminus \{j(i)\}}=0,
\end{equation}
since the middle term is
\begin{eqnarray*}
 \sum_{i=1}^{r} (-1)^{i+1} a^{q^{r-2}}_{j(i)}  \left( \sum_{k=1}^{i-1} (-1)^{k+1}a^{q^{r-2}}_{j(k)} h^q_{J \setminus \{j(i),j(k)\}}\right. \hspace{4.1cm}\\
 \hspace{-1cm}  +\left. \sum_{k=i+1}^{r} (-1)^{k} a^{q^{r-2}}_{j(k)} h^q_{J \setminus
  \{j(i),j(k)\}}\right)=\hspace{1cm}\\
\sum_{1 \leq k<i \leq r}\!\!\!\!\!(-1)^{i+k}
 a^{q^{r-2}}_{j(i)} a^{q^{r-2}}_{j(k)} h^q_{{J \setminus \{j(i),j(k)\}}}
 +\!\!\!\!\sum_{1 \leq i<k \leq r}\!\!\!\!\!(-1)^{i+k+1}
a^{q^{r-2}}_{j(i)} a^{q^{r-2}}_{j(k)} h^q_{{J \setminus \{j(i),j(k)\}}},
\end{eqnarray*}
 which equals zero. When $r=n+1$, $u_{J}\in Z^{0}(G,\kv)$. It follows that $u_J(\iota) \in \kvg$, which implies
  $u_{J}(\iota)^{q}\in A$ ($\iota\in G$ is the neutral element). Therefore, for each $1 \leq i \leq n+2$, we have
  $u^q_{\{1,2,\ldots n+2\}\setminus \{i\}}(\iota) \in A$. The second equality in (\ref{uJdnr2}) is also valid
for  $J=\{1,2,\ldots,n+2\}$ ($r=n+2$), that is, 
\[\sum_{i=1}^{n+2}(-1)^{i+1} a^{q^n}_i u^q_{\{1,2,\ldots, n+2\}\setminus \{i\}}(\iota) = 0.\]
As $a^{q^n}_1, \ldots, a^{q^n}_{n+2}$ is $A$-regular,  $u^q_{\{1, \ldots, n+1\}}(\iota) \in (a^{q^n}_1,\ldots, a^{q^n}_{n+1})A$. Thus there exist $f_1,f_2, \ldots, f_{n+1}\in \Hom_{{\mathbb Z}G}(X_{0},A)$ (in particular, $f_{i}\partial_{1}=0$) such that $u^q_{\{1,\ldots,n+1\}}= \sum_{i=1}^{n+1}a^{q^n}_i f_i$. Substituting (\ref{DefNOfuJ})
for $u_{\{1, \ldots, n+1\}}$ yields 
$$\sum_{i=1}^{n+1} (a_i^{q^n})((-1)^{i+1} h^q_{\{1, \ldots, n+1\} \setminus
  \{i\}}-f_i) = 0.$$ 
Since $a^{q^n}_1, \ldots, a^{q^n}_{n+1}$ is a phsop in
$A$, it is also a phsop in $\kk[V]$, and so $((-1)^n h^q_{\{1,\ldots, n\}}-f_{n+1})(\iota) \in (a^{q^n}_1,a^{q^n}_2, \ldots,
a^{q^n}_n)\kk[V]$. It follows that 
\[h^q_{\{1, \ldots, n\}}- (-1)^n
f_{n+1} = \sum_{i=1}^n a^{q^n}_i l_i\] 
for some $l_1, l_{2},\ldots, l_{n} \in\Hom_{\mathbb{Z}G}(X_{0},\kk [V]) \cong \kv$. We next apply $\partial_1$
  to this expression. For $n=1$, we have
  $a_{1}^{q}u_{\{1\}}^{q}=a_{1}^{q}l_{1}\partial_{1}$, or alternatively
  $u_{\{1\}}^{q}=l_{1}\partial_{1}\in B^{1}(G,\kv)$, a contradiction to $g^{q}\ne
  0$. Now assume that $n\ge 2$.  Applying $\partial_1$ leads to\begin{equation}\label{reqn}
u_{\{1, \ldots , n\}}^{q^2} = \sum_{i=1}^n a^{q^n}_i l_i \partial_1.\end{equation}
For $2 \leq r \leq n$, we prove by reverse induction that there exist elements  $l_1,l_2,\ldots,l_r\in\Hom_{\mathbb{Z}G}(X_{n-r},\kk [V])$ such that 
\begin{equation}\label{uqnr}
u^{q^{n+2-r}}_{\{1, \ldots, r\}} = \sum_{i=1}^r a_i^{q^n}
l_i\partial_{n-r+1}.\end{equation} 
The case $r=n$ is covered by \eqref{reqn}. Suppose $u^{q^{n+1-r}}_{\{1,
  \ldots, r+1\}} = \sum_{i=1}^{r+1} a^{q^n}_i l'_i \partial_{n-r}$ for some $l'_1,l'_2, \ldots l'_{r+1}\in \Hom_{\mathbb{Z}G}(X_{n-r-1},\kk[V])$. Using (\ref{DefNOfuJ}), we obtain
$$\sum_{i=1}^{r+1} (-1)^{i+1} a^{q^{n}}_i h^{q^{n+1-r}}_{\{1,\ldots, r+1\} \setminus \{i\}} = \sum_{i=1}^{r+1} a^{q^n}_i l'_i \partial_{n-r},$$ and rearranging yields
$$\sum_{i=1}^{r+1} a^{q^{n}}_i ((-1)^{i+1}  h^{q^{n+1-r}}_{\{1,\ldots, r+1\} \setminus \{i\}} - l'_i \partial_{n-r})=0.$$ 
Since $a^{q^n}_1, \ldots, a^{q^n}_{r+1}$ is a phsop for $\kv$, we have
$$ h^{q^{n+1-r}}_{\{1, \ldots, r\}} - (-1)^r l'_{r+1} \partial_{n-r} =
\sum_{i=1}^r a_i^{q^n} l_i$$ 
for some $l_1,l_2, \ldots, l_{r}\in\Hom_{\mathbb{Z}G}(X_{n-r},\kk [V])$. Here we have used that $X_{n-r}$ is a free
${\mathbb Z}G$-module. Applying $\partial_{n-r+1}$ to this expression gives us
$$h^{q^{n+1-r}}_{\{1, \ldots, r\}} \partial_{n-r+1} = \sum_{i=1}^r a_i^{q^n} l_i
\partial_{n-r+1},$$ 
which implies \eqref{uqnr}, as required. 

When $r=2$, Equation~\eqref{uqnr} reads
$u^{q^n}_{\{1,2\}} = a^{q^n}_1 l_1 \partial_{n-1} + a^{q^n}_2 l_2
\partial_{n-1}$, where $l_1,l_2 \in
\Hom_{\mathbb{Z}G}(X_{n-2},\kk[V])$. Substituting (\ref{DefNOfuJ}) for $u_{\{1,2\}}$, we obtain $a^{q^n}_1 ( h^{q^n}_{\{2\}} - l_1 \partial_{n-1}) = a^{q^n}_2
(h_{\{1\}}^{q^n} + l_2 \partial_{n-1})$. For the same reasons as before, we have
$h^{q^n}_{\{1\}} + l_2 \partial_{n-1} = a_1^{q^n} l$ for some $l \in
\Hom_{\mathbb{Z}G}(X_{n-1},\kk[V])$. Applying $\partial_n$ gives us $a_1^{q^n}
u_{\{1\}}^{q^n} = a_1^{q^n} l \partial_n$. Hence, $u_{\{1\}}^{q^n} = l
\partial_n\in B^{n}(G,\kv)$, a contradiction to $g^{q^{n}}\ne 0$. Therefore, $\depth_{A}(I)\le n+1$.

Finally, if $c_{1},\ldots,c_{k}\in \Ann_{\kvg}(g)$ forms a phsop in $\kvg$,
then $c_{1}^{q},\ldots,c_{k}^{q}\in I$ forms a phsop in $A$, and so
$I$ has height at least $k$. The graded analogue of \cite[Exercise
1.2.23]{BrunsHerzog} implies the Cohen-Macaulay defect of $A$ is at least $\height_{A}(I)-\depth_{A}(I)\ge k-n-1$.
\end{proof}

%----------------------------------------------------------------

\begin{Lemma}\label{trivcoeff}
If $g\in H^{n}(G,\kk)$ is nonzero, then $g^{p^{m}}$ is nonzero for all
$m\geq 0$. 
\begin{proof}
 For
$n=1$, this is clear since elements  of $H^{1}(G,\kk)$ are group homomorphisms
$G\rightarrow(\kk,+)$. For arbitrary $n$, by the Universal Coefficient
Theorem \cite[page
30]{Evens}, $H^{n}(G,\kk)\cong H^{n}(G,\FF_{p})\otimes_{\FF_{p}}\kk$. We have $g=\sum_{i\in I} g_{i}\otimes \lambda_{i}$ for some $g_{i}\in
H^{n}(G,\FF_{p})$ and $\lambda_{i}\in \kk$. Without loss of generality, we can assume that
the set $\{\lambda_{i}: i\in I\}$ is $\FF_{p}$-linearly independent and $g_{i}$ is
nonzero for all $i$. Then $g^{p^{m}}=\sum_{i\in I} g_{i}^{p^{m}}\otimes
\lambda_{i}^{p^{m}}=\sum_{i\in I} g_{i}\otimes \lambda_{i}^{p^{m}}$, since the
$m$-fold Frobenius homomorphism induces the identity map on $\FF_{p}$, and
thus also on $H^{1}(G,\FF_{p})$. Therefore,  as $\{\lambda_{i}^{p^{m}}: i
\in I\}$ is still $\FF_{p}$-linearly independent, $g^{p^{m}}$ is still
nonzero.
\end{proof}
\end{Lemma}

\begin{rem}
 For $n>1$, Theorem \ref{MainTheorem} is new even in the case $A=\kk[V]^G$. 
\end{rem}

\begin{eg}
Let $G\subseteq (\kk,+)$ be a finite nontrivial subgroup. Consider
the threefold sum of the $2$-dimensional representation $V$ of $G$ over $\kk$ given by \!\!\!\!
\renewcommand{\arraystretch}{0.7}
\renewcommand{\arraycolsep}{2pt}
$\sigma\mapsto
{\tiny
{\left(\begin{array}{cc}1&0\\-\sigma
        &1\end{array}\right)}}$. Write $\kk[V^{\oplus
  3}]=\kk[x_{1},y_{1},x_{2},y_{2},x_{3},y_{3}]$ and $\rho$ for the induced $G$-action. The map $g: G\rightarrow
\kk,\,\, \sigma\mapsto \sigma$ yields a nonzero element in $H^{1}(G,\kk)$. For
all $\sigma\in G$, we have $x_{i}g_{\sigma}=(\rho(\sigma)-1)y_{i}$, that is, $x_{i}g$ is trivial. Since $x_{1},x_{2},x_{3}$ form a phsop in
$\kk[V^{\oplus 3}]$, Theorem~\ref{MainTheorem} implies that no graded
geometric separating algebra in $\kk[V^{\oplus 3}]^{G}$ is Cohen-Macaulay.\done
\end{eg}

The following example shows that Theorem~\ref{MainTheorem} applies only to
graded \emph{geometric} separating algebras:

\begin{eg}\label{C4perm}
Let $V$ be the permutation representation of the cyclic group
$C_{4}=\langle \sigma\rangle$ of order $4$ over the field $\FF_2$. Consider the $C_4$-invariants
 $c_{1}:=x_{1}+x_{2}+x_{3}+x_{4}$, $c_{2}:=x_{1}x_{3}+x_{2}x_{4}$,
$c_{3}:=x_{1}x_{2}+x_{2}x_{3}+x_{3}x_{4}+x_{1}x_{4}$, and $c_4:=x_1x_2x_3x_4$. 
The action of $C_4$ on $V$ partitions its $16$ elements into $6$ orbits, which one can check are separated by $c_1$, $c_2$, $c_3$, and $c_4$. As
$c_{1},c_{2},c_{3},c_{4}$ form a hsop in $\FF_2[V]$, the subalgebra $\FF_2[c_1,c_2,c_3,c_4]$ is a polynomial graded (non geometric) separating algebra.  In particular, it is Cohen-Macaulay.

On the other hand, if $g: C_4\rightarrow \FF_2$ is
the nontrivial cocycle given by  $\sigma^{i}\mapsto i \mod 2$, then
$c_{1}g_{\sigma^{i}}=(\sigma^{i}-1)(x_{1}+x_{3})$,
$c_{2}g_{\sigma^{i}}=(\sigma^{i}-1)(x_{1}x_{3})$, and
$c_{3}g_{\sigma^{i}}=(\sigma^{i}-1)(x_{1}x_{4}+x_{2}x_{3})$, that is, $c_1,c_2,c_3\in\Ann_{\overline{\FF}_{2}[\overline{V}]^{C_4}}(g)$. Since
$c_{1},c_{2},c_{3},c_{4}$ form a hsop in $\FF_2[V]$, by Theorem \ref{MainTheorem}, no graded geometric separating algebra in $\FF_2[V]^{C_4}$ is  Cohen-Macaulay.\done
\end{eg}

Using Theorem~\ref{MainTheorem}, we can quickly generalize several results
which were consequences of its analogue \cite[Corollary~1.6]{KemperOnCM}. A
first consequence is the following generalization of \cite[Corollary
21]{CampEtAl}:

\begin{Corollary}\label{3copiesTheorem}
Assume $G$ contains a normal subgroup $N$ of index $p$ (for example, $G$ is a $p$-group). Then for any faithful representation $V$, every graded geometric separating algebra in  $\kk[V^{\oplus k}]^{G}$ has Cohen-Macaulay defect at least $k-2$.
\begin{proof} Since $G/N\cong ({\mathbb F}_{p},+)$, there is a nonzero element in $H^{1}(G,\kk)$. As $V$ is faithful, the fixed subspaces of nonidentity elements in $G$ have codimension at least $k$ in $V^{\oplus k}$. The result follows from Lemma \ref{GenericAnnihilators}. \end{proof}\end{Corollary}

\begin{Lemma}\label{GenericAnnihilators}
Suppose $V$ is faithful. If the fixed subspace of every element of order $p$ in $G$ has codimension at least $k$ in $V$, then for any homogeneous $g\in
  H^{n}(G,\kvbar)$, the ideal $\Ann_{\kvgbar}(g)$ has height at least~$k$.

Therefore, if there exists a homogeneous $g\in
  H^{n}(G,\kvbar)$ satisfying the hypotheses of Theorem \ref{MainTheorem}, then every graded geometric separating algebra in $\kvg$ has Cohen-Macaulay defect at least $k-n-1$. 
\begin{proof} We may assume $\kk=\kbar$. By Kemper \cite[Lemma 2.1]{KemperOnCM} there exist (in general, non-homogeneous) elements $a_{1},\ldots,a_{k}\in \kvgplus$ such that $a_{i}H^{n}(G,\kv)=0$, for all $i$, and $\height (a_{1},\ldots,a_{k})=k$. It follows that $\Ann_{\kvg}(g)$ has height at least $k$.\end{proof}\end{Lemma}

\begin{proof}[Proof of Theorem \ref{kcopiesTheorem}]
By \cite[Theorem~4.1.3]{BensonCohomology2}, there is a number $r$ such that
$H^{r}(G,\kk)\ne 0$. Thus by Lemma~\ref{trivcoeff}, for any $k$, there is a
minimal number $n\le r$ such that the hypotheses of Theorem \ref{MainTheorem}
are satisfied for $V^{\oplus k}$. The same argument as in Corollary
\ref{3copiesTheorem} shows that $A$ has Cohen-Macaulay defect at least
$k-n-1\ge k-r-1$.
\end{proof}

Next, we generalize three results of \cite{KemperOnCM}. Note that since, for us, elements acting trivially are bireflections, we do not need to assume that $V$ is faithful.

\begin{Corollary}
If $G$ has 
a normal subgroup $N$ of index $p$ which contains all elements acting as bireflections on $V$,
then no graded geometric separating algebra in $\kvg$ is Cohen-Macaulay.
\begin{proof} The proof of \cite[Theorem~3.6]{KemperOnCM} shows  that the
hypotheses of Theorem \ref{MainTheorem} are fulfilled with $n=1$ and $k=3$. 
\end{proof}\end{Corollary}

\begin{proof}[Proof of Theorem~\ref{thm1}]
For $p$-groups, if the elements acting as bireflections generate a proper
subgroup, then this subgroup lies in a normal subgroup of index $p$.\end{proof}

\begin{Prop}\label{bireflectionCriterion}
Suppose $G$ has a normal subgroup $N$ with factor group an
elementary abelian $p$-group. If there exists $\sigma\in G\setminus
N$ whose fixed space in $V$ is not contained in the fixed space of any bireflection in $G\setminus N$, then no graded geometric separating algebra in $\kvg$ is Cohen-Macaulay.
\begin{proof} Without loss of generality, assume $\kk=\kbar$. As $G/N$ is an elementary abelian $p$-group, there is a $g\in H^{1}(G,\kk)$ with kernel $N$. The proof of \cite[Theorem~3.9]{KemperOnCM} provides a phsop $a_{1},a_{2},a_{3}$ in $\kvg$\!\!, and a homogeneous invariant $h\notin\sqrt{\Ann_{\kvg}(g)}$ such that $ha_{i}\in\sqrt{\Ann_{\kvg}(g)}$, for $i=1,2,3$. Hence, for some $k\geq 0$, the invariants  $a_{i}^{k}$ annihilate $g':=h^{k}g\in H^{1}(G,\kv)$.
Thus, the annihilator ideal $\Ann_{\kvg}(g')$  has height at least $3$. By Theorem~\ref{MainTheorem}, it now suffices to show that $g'^{p^{m}}$ is nonzero for all $m\geq 0$. By \cite[Proposition~3.5]{KemperOnCM}, we have $\sqrt{\Ann_{\kvg}(g^{p^{m}})}=\sqrt{\Ann_{\kvg}(g)}$. Therefore, if $g'^{p^{m}}$ is zero, then $h^{kp^m}$ annihilates $g^{p^m}$, and so $h\in\sqrt{\Ann_{\kvg}(g)}$, a contradiction.
\end{proof}\end{Prop}

\begin{eg}\label{refnEg}
Let $\kk$ be a finite field. For $m\ge 3$, set $V=\kk^{2m+1}$, and consider the group $G\leq\GL(V)$ formed by the $({2m+1})\times ({2m+1})$ matrices of the form
\[\renewcommand{\arraystretch}{0.7}
\renewcommand{\arraycolsep}{2pt}
\left(\begin{array}{c|c}
I_{m+1} & \mathbf{0}\\
\hline
\begin{array}{cccc}
\alpha_{0}&&&\alpha_{m}\\
&\ddots&&\vdots\\
&&\alpha_{m-1}&\alpha_{m}
\end{array}
&
I_{m}
\end{array}\right),
\]
where $\alpha_{0},\ldots,\alpha_{m}\in \kk$, and  $I_{m}$ denotes the
$m\times m$ identity matrix. The group $G$ is a $p$-group, and is generated by
reflections, that is, by elements whose fixed space has codimension at most 1 in $V$. Example~3.10 in \cite{KemperOnCM} shows that the hypotheses of
Proposition~\ref{bireflectionCriterion} are satisfied, with $N$ the subgroup
formed by the elements such that $\alpha_m=0$, and $\sigma$ the element such that $\alpha_i=1$ for all $i$. Hence, no graded geometric separating algebra in $\kvg$ is Cohen-Macaulay.\done
\end{eg}

We end this section with a generalization of \cite[Theorem~2.7]{KemperOnCM}.

\begin{Theorem}\label{Vreg}
Let  $V_{reg}$ be the regular representation of $G$ over $\kk$. If $p$ divides $|G|$, then every graded geometric separating algebra in $\kk[V_{reg}]^G$ has Cohen-Macaulay defect at least $|G|\frac{p-1}{p}-2$. For $|G|\ge 5$, this number is at least one.
\begin{proof}
By \cite[Lemma 2.6]{KemperOnCM}, there is a nonzero $g\in H^1(G,\kk[V_{reg}])$. Lemma~\ref{RegulRepCocycle} implies all powers $g^{p^m}$ are also nonzero. Since the fixed subspaces of elements of $G$ of order $p$ have codimension $|G|(p-1)/p$, the hypotheses of Lemma~\ref{GenericAnnihilators} are satisfied with $k=|G|(p-1)/p$.
\end{proof}
\end{Theorem}

The regular representation of $C_4$ was studied in Example \ref{C4perm}. For
$G=C_2,C_3$, and $C_2\times C_2$, the invariant ring $\kk[V_{reg}]^{G}$ is
Cohen-Macaulay \cite[Theorem~2.7]{KemperOnCM}. Thus, these are the only groups
such that there exists a Cohen-Macaulay graded geometric separating algebra in
$\kk[V_{reg}]^{G}$.
\begin{Lemma}\label{RegulRepCocycle}
Let $V$ be a permutation representation of $G$. If $g$ in $H^{1}(G,\kv)$ is nonzero, then $g^{p^m}$ is nonzero for all $m\ge 0$. If, in addition, $V$ is faithful, then every graded geometric separating algebra in $\kk[V^{\oplus k}]^G$ has Cohen-Macaulay defect at least $k-2$.
\begin{proof}
As $V$ is a permutation representation, there is a set of monomials
$M\subseteq\kv$ such that $\kv=\bigoplus_{h\in M} \langle Gh \rangle$. Thus,
if $g\in H^{1}(G,\kv)$, there is a (finite) decomposition $g=\sum_{h\in M}
g_h$, where each $g_h$ is in $H^{1}(G,\langle Gh \rangle)$. For  $m\ge 0$, we
have the decomposition $g^{p^{m}}=\sum_{h\in M} g_h^{p^m}$\!\!, where $g_h^{p^m}$
is in $H^1(G,\langle Gh' \rangle)$, and $h'$ is the unique element of
$Gh^{p^m}\cap M$. If $g$ is nonzero, then $g_h$ is nonzero for some $h$. As
$\langle Gh \rangle$ and $\langle Gh^{p^m}\rangle$ are isomorphic permutation
representations of $G$, the element $g_h^{p^m}$ is also nonzero. The additional statement follows by Lemma \ref{GenericAnnihilators}.
\end{proof}
\end{Lemma}

%%%%%%%%%%%%%%%%%%%%%%%%%%%%%%%%%%%%%%%%%%%%%%%%%%%%%%%%%%%%%%%%%%%%%%%%%%%%%%%%%%%%%%%%%%%%%%%%%%%%%%%%%%%%%%%%
%%%%%  Section:      A4              %%%%%%%%%%%%%%%%%%%%%%%%%%%%%%%%%%%%%%%%%%%%%%%%%%%%%%%%%%%%%%%%%%%%%%%%%%%%%%%%%%%%%%%%%%%%%%%%%%%%%%%%%%%%%%%%%%%%%%%%%%%%%%%%%%%%%%%%%%%%%%%%%%%%%%%%%%%%%%%%%%%%%%%%%%%%%%%%%%%%%%%%%%%%%%%%%%%%%%%%%%%%%%%%%%

\section{The Alternating Group $A_4$}\label{A4}

In this section we concentrate on representations of the alternating group
$A_4$ over an algebraically closed field $\kk$ of characteristic $2$. The
group $A_4$ is the smallest modular group $G$ for which the cohomology
$H^1(G,\kk)$ is trivial. In order to apply  Theorem \ref{MainTheorem}, we must look outside the direct summand $H^1(G,\kk)$ of $H^1(G,\kk[V])$. 

Pick $\chi\in A_4$, an element of order 2, and $\tau$, a 3-cycle, so that $A_4$
is generated by $\chi$ and $\tau$. The unique Sylow $2$-subgroup $P$ in $A_4$
is generated by $\chi$ and $\tau^{-1}\chi\tau$. Let $\omega$ be a fixed
primitive third root of unity in $\kk$. Define $^{\omega}\kk$ and
$^{\omega^2}\kk$ to be the one-dimensional  representations on which $P$ acts
trivially, and $\tau$ acts via multiplication by $\omega$ and $\omega^{2}$,
respectively. Formally, $^1\kk$ denotes the trivial representation.

\begin{Lemma}\label{ktwiddles}
For $i=1,2$, the element of $H^1(A_4, ^{\omega^i}\hspace{-1mm}\kk)$ given by
the cocycle $\chi \mapsto \omega^{2i}$, $\tau\mapsto 0$ is nonzero.
\begin{proof}
Assume $i=1$.  Let $V:= \langle v_1,v_2 \rangle$ be the representation given by \renewcommand{\arraystretch}{0.7}
\renewcommand{\arraycolsep}{2pt}
\[
\chi\mapsto\left(
\begin{array}{cc}
1 & \omega^2 \\
0 & 1 
\end{array}
 \right),~\mathrm{and}\,\,\,
\tau\mapsto\left(
\begin{array}{cc}
 \omega & 0 \\
0 & 1 
\end{array}
\right). 
\]
  This is an indecomposable representation \cite[Theorem 7.0.3]{ElmerThesis},
  with a submodule $W:= \langle v_1 \rangle$ isomorphic to
  $^{\omega}\kk$. Furthermore, $A_4$ acts trivially on the quotient $V/W$, and so
  we have a nonsplit exact sequence
$$
0 \rightarrow ^{\omega} \hspace{-1mm} \kk \rightarrow V \rightarrow \kk \rightarrow 0.
$$
By \cite[section 2]{KemperLinRed}, it follows that $H^1(A_4, ^{\omega}\kk)
\neq 0$. In particular, the cocycle  given by $\chi \mapsto \omega^{2}$,
$\tau\mapsto 0$ is nontrivial. Similarly, for $i=2$, the cocycle in $Z^1(A_4, ^{\omega^{2}}\kk)$ given by  $\chi \mapsto \omega$, $\tau \mapsto 0$ is nontrivial.\end{proof}\end{Lemma}

\begin{Corollary}
If $V:=  ^{\omega^i} \hspace{-1.5mm} \kk \oplus W^{\oplus k}$, where $W$ is a faithful representation and $i=1,2$, then every graded geometric separating algebra in $\kv^{A_4}$ has  Cohen-Macaulay defect at least $k-2$.
\begin{proof}
Without loss of generality, $i=1$. Since $V^* \cong ^{\omega^2} \hspace{-1mm}
\kk \oplus (W^*)^{\oplus k}$, there is  a direct summand of $\kv$ isomorphic
to $S( ^{\omega^2}\kk)$. Since $S^{2^m}(^{\omega^2}\kk)$ is isomorphic to $ ^{\omega^2}\kk$ or $ ^{\omega}\kk$,
Lemma~\ref{ktwiddles} implies that the cohomology class $g_{m} \in H^1(A_4,S^{2^m}(
^{\omega^2}\kk))$, given by the cocyle $\chi \mapsto \omega^{2^{m}}$, $\tau
\mapsto 0$, is nonzero for all  $m$. As $g_{m}=g_{0}^{2^{m}}$, the
result follows by Lemma~\ref{GenericAnnihilators}. 
\end{proof}\end{Corollary}

Exploiting the classification of the finite dimensional representations of $A_4$~\cite{Conlon} (with the notation of \cite[Chapter~7]{ElmerThesis}), we obtain a much stronger result:

\begin{Theorem} Suppose $V$ is a faithful, indecomposable finite dimensional representation of $A_4$. If $\kk[V]^{A_4}$ is non Cohen-Macaulay, then no graded geometric separating algebra in $\kk[V]^{A_4}$ is Cohen-Macaulay.
\begin{rem}
The indecomposable representations $V$  of ${A_4}$ such that $\kk[V]^{A_4}$ is Cohen-Macaulay are listed in \cite[Corollary~5.2.16]{ElmerThesis}. In particular, when $\dim_\kk(V) \geq 7$, $\kk[V]^{A_4}$ is non Cohen-Macaulay. 
\end{rem}
\begin{proof}
Whenever $\kk[V]^{A_4}$ is non Cohen-Macaulay, the fixed space of any
nonidentity element in $P$ has codimension at least $3$ in V (see the classification). Thus, by Lemma~\ref{GenericAnnihilators}, it suffices to find $g \in H^1({A_4},\kk[V])$ such that $g^{2^m}$ is nonzero for all $m\geq 0$. We use the classification to separate our argument into two cases. 

First, we suppose $V$ is of the form $W_s(\omega^e)$,
where $s\in\ZZ$ ($s$ may be negative only if it is odd), and $0 \leq e
\leq 2$. The dimension of $V$ is $n:= |s|$. Let $\{v_1,v_2, \ldots,
v_n\}$ be a basis of $V$ such that $A_4$ acts via  
the matrices given in \cite[Theorem~7.0.3]{ElmerThesis}. In particular, the
action of ${A_4}$ is upper triangular, and $\tau$ acts diagonally via
multiplication by third roots of unity. If $\kk[V]^{A_4}$ is non
Cohen-Macaulay, the subset $\{v_1,v_2, \ldots, v_l\}$ is contained in $V^P$,
for some $l \geq 2$. Moreover for some $r \leq l$, we have $\tau v_r =
\omega^{-j} v_r$, where $j=1$ or 2. Let  $\{x_1,x_2, \ldots, x_n\}$ be the
dual basis, so that $\kk[V] = \kk[x_1, x_2, \ldots x_n]$.  Define $u:=
\prod_{\sigma \in P}(\sigma x_r)$. Then 
$$\tau u = \prod_{\sigma \in P}(\tau\sigma \tau^{-1} (\tau x_r)) = \prod_{\sigma \in P}\omega^j (\sigma x_r) =\omega^j u .$$ 
Since $v_r \in V^P$, the set $\{x_1,x_2, \ldots, x_{r-1},x_{r+1}, \ldots, x_n\}$ spans a $\kk
P$-submodule of $V^*$, thus $\kk[V] = \kk[x_1, x_2, \ldots,
x_{r-1},x_{r+1}, \ldots, x_n][x_r]$ as $\kk P$-modules. Therefore,  by
\cite[Proof of Lemma~1.3.2]{ElmerThesis}, if $M$ is a $\kk P$-module direct
summand of $S(V^*)$, so is $u \cdot M$, and by induction so is $u^i \cdot M$
for any integer $i$. In particular, as $1 \in S^0(V^*)$, we have a sequence $\langle 1 \rangle, \langle u \rangle, \langle u^2 \rangle, \ldots $ of $\kk P$-direct summands. Each $\langle u^i \rangle$ is a $\kk A_4$-direct summand, and $\langle u^i \rangle \cong ^{\omega^{ij}}\hspace{-1mm} \kk$.
The element $g_{m} \in H^1({A_4},\langle u^{2^{m}} \rangle)$ given by the
cocyle $\chi \mapsto \omega^{-2^m j} u^{2^m}$, $\tau \mapsto 0$, is nonzero by
Lemma \ref{ktwiddles}. As $g_{m}=g_{0}^{2^{m}}$, we are done.

Second, suppose $V$ takes the form $\overline{W}_{6d, \lambda}$, where $\lambda \in \kk \cup \{\infty\}$, not a
third root of unity, and $d \geq 1$, an integer, are such that if $d=1$, then $\lambda$ is not $0$, $1$, or $\infty$ (see \cite[Theorem~7.0.3]{ElmerThesis} for
notation). As a $\kk P$-module, $V^* \cong V_{2d,\alpha} \oplus
V_{2d, \beta} \oplus V_{2d, \gamma}$ \cite[Section~4.3]{BensonCohomology1},
where $\alpha$, $\beta$, and $\gamma$ are elements of $\kk \cup \{\infty\}$
depending on $\lambda$. The element $\tau$ acts on $V$ by permuting the three
$\kk P$-module summands. As a $\kk P$-module, $S(V^*) \cong
S(V_{2d,\alpha}) \otimes S(V_{2d, \beta}) \otimes S(V_{2d, \gamma})$, and
therefore, $S(V^{*})$ has direct summands isomorphic to $S(V_{2d,\alpha})$, $S(V_{2d, \beta})$, and $S(V_{2d, \gamma})$.
These summands are permuted by the action of $\tau$. Let $\{x_1, \ldots, x_{6d}\}$ be a
basis of $V^*$ such that the action of $P$ is given in block form; that is, $x_1, \ldots, x_{2d}$ is a basis for the summand
isomorphic to $V_{2d, \alpha}$, and so on. Define $u_{\alpha}:= \prod_{\sigma \in
  P}\sigma  x_1$, and  $u_{\beta}:= \prod_{\sigma \in P}\sigma x_{2d+1}$, and
$u_{\gamma}:= \prod_{\sigma \in P}\sigma  x_{4d+1}$. By the same argument as
before, for each integer $i$, $\langle u^i_{\alpha} \rangle$ is a direct
summand (with trivial $P$-action) of the $\kk P$-module $S(V_{2d,\alpha})$,
and hence of the $\kk P$-module $S(V^*)$. The analogue holds for $u_{\beta}$
and $u_{\gamma}$. Furthermore, we have 
$$\tau  u^i_{\alpha} = \tau (\prod_{\sigma \in P} \sigma  x_1)^i =
(\prod_{\sigma \in P} (\tau \sigma \tau^{-1})   x_{2d+1})^i = u^i_{\beta}.$$ 
Similarly,  $\tau u_{\beta}^{i}=u_{\gamma}^{i}$ and $\tau
u_{\gamma}^{i}=u_{\alpha}^{i}$. Thus, $U^i:=\langle u^i_{\alpha}, u^i_{\beta},
u^i_{\gamma} \rangle$ is a $\kk {A_4}$-direct summand of $S(V^*)$ on which $P$
acts trivially. Since $[{A_4}:P]$ is odd, $\tau$ acts diagonally with respect to some basis of $U^{i}$. A short calculation shows that $u:= u_{\gamma}+ \omega u_{\beta}+ \omega^2 u_{\alpha}$ spans a direct summand of $U^{1}$ on which $\tau$ acts via multiplication by $\omega$, in other words, $\langle u \rangle \cong ^{\omega}\hspace{-1mm}\kk$. For each $m$, we have that $u^{2^m} = u^{2^m}_{\gamma}+\omega^{2^m} u^{2^m}_{\beta}+\omega^{2^{m+1}}u^{2^m}_{\alpha}$ spans a direct summand of $U^{2^m}$ isomorphic to $^{\omega^{2^m}} \kk$.

 The element $g_{m} \in H^1({A_4},\langle u^{2^m} \rangle)$ given by the
 cocyle $\chi \mapsto \omega^{2^{m+1}} u^{2^m}$, $\tau \mapsto 0$ is nonzero
 for each $m \geq 0$. As $g_{m}=g_{0}^{2^{m}}$, the result now follows by Lemma \ref{GenericAnnihilators}.
\end{proof}\end{Theorem}

%%%%%%%%%%%%%%%%%%%%%%%%%%%%%%%%%%%%%%%%%%%%%%%%%%%%%%%%%%%%%%%%%%%%%%%%%%%%%%%%%%%%%%%%%%%%%%%%%%%%%%%%%%%%%%%%%%%%%%%%%%%%%%%%%%%%%%%%%%%%%%%%%%%%%%%%%%%%%%%%%%%
%%%%      SECTION:     Concluding remarks                       %%%%%%%%%%%%%%%%%%%%%%%%%%%%%%%%%%%%%%%%%%%%%%%%%%%%%%%%%%%%%%%%%%%%%%%%%%%%%%%%%%%%%%%%%%%%%%%%%%%%%%%%%%%%%%%%%%%%%%%%%%%%%%%%%%%%%%%%%%%%%%%%%%%%%%%%%%%%%%%%%%%%%%%%%%%%%%%%%%%%%%%%%%%%%%%%%%%

\section{Concluding Remarks}\label{remarks}

Our results have shown that in many of the situations in which the ring of
invariants is non Cohen-Macaulay, the same holds for any graded geometric
separating algebra. Hence, one might wonder if the existence of a Cohen-Macaulay
graded geometric separating algebra implies that the ring of invariants itself
is Cohen-Macaulay. The following example shows this is not true.

\begin{eg}
Let $G=C_2 \times C_2=\langle \sigma,\tau\rangle$ be the Klein four group, and let $\kk$ be a field of characteristic two. Consider the 5-dimensional representation of $G$ given by 
\renewcommand{\arraystretch}{0.7}
\renewcommand{\arraycolsep}{2pt}
\[\sigma\mapsto {\tiny\left(\begin{array}{ccccc} 1&0&1&0&0\\0&1&0&1&0\\0&0&1&0&0\\0&0&0&1&0\\0&0&0&0&1\end{array}\right)},\ \tau\mapsto{\tiny\left(\begin{array}{ccccc}1&0&0&1&0\\0&1&0&0&1\\0&0&1&0&0\\0&0&0&1&0\\0&0&0&0&1\end{array}\right)}.\]
By \cite[Theorem~7]{ElmerFleischmann}, the ring of invariants
is not Cohen-Macaulay. In fact, \cite{ElmerFleischmann} shows that there
exists a nonzero cohomology class $g \in H^1(G,\kk[V])$ whose restriction to
each proper subgroup of $G$ is zero. By \cite[Lemma~2.4,
Proposition~2.5]{ElmerAssocPrimes}, $\sqrt{\Ann_{\kvg}(g)} = \frak{I}(V^G)
\cap \kvg$, where  $\frak{I}(V^{G})$ denotes the ideal of polynomials $f \in \kv$ vanishing on $V^{G}$. Thus, the height of $\Ann_{\kk[V]^G}(g)$ is $\codim(V^G) = 3$ while, by \cite[Corollary~1.6]{KemperOnCM}, its depth is only two.

Using MAGMA \cite{Magma} and the methods of \cite[Section 2]{KemperCompRed}, one can verify that \[ \begin{array}{l} \{ a_{1}:= x_3, 
    a_{2}:= x_4, 
    a_{3}:= x_5,\\
   \ a_{4}:= x_1^4 + x_1^2 x_3^2 + x_1^2 x_3 x_4 + x_1 x_3^2 x_4 + x_1 x_3 x_4^2 + x_1 x_3 x_4 x_5 + \\
    \hspace{1cm}    x_1 x_4^3 + x_2^2 x_3^2 + x_2 x_3^2 x_5 + x_2 x_3 x_4^2,\\
 \ a_{5}:= x_2^4 + x_2^2 x_4^2 + x_2^2 x_4 x_5 + x_2^2 x_5^2 + x_2 x_4^2 x_5 + x_2 x_4 x_5^2,\\
   \ a_{6}:= x_1^2 x_4^2 + x_1 x_3 x_4 x_5 + x_1 x_4^3 + x_2^2 x_3^2 + x_2 x_3^2 x_5 + x_2 x_3 x_4^2,\\
   \ a_{7}:= x_1 x_4^2 x_5 + x_1 x_4 x_5^2 + x_2^2 x_3 x_5 + x_2^2 x_4^2 + x_2 x_3 x_5^2 + x_2 x_4^3\}  \end{array}\]
forms a geometric separating set. Furthermore, $\{a_{1},a_{2},a_{3},a_{4},a_{5}\}$ is a hsop for the geometric separating
    algebra $A: =\kk[a_{1},a_{2},\ldots,a_{7}]$. As a module over
    $\kk[a_{1},a_{2},a_{3},a_{4},a_{5}]$, $A$ is freely generated by
    $\{1,a_{6},a_{7},a_{6}a_{7}\}$. Therefore, $A$ is Cohen-Macaulay. As $G$ is a
    $p$-group, by Theorem
    \ref{thm1}, $G$ must be a bireflection group, which is indeed the case. The Hilbert Series of $A$ is
    $H(A,t)=\frac{1+2t^{4}+t^{8}}{(1-t)^{3}(1-t^{4})^{2}}$. Since $H(A,1/t)=(-1)^{5}t^{3}H(A,t)$, we even have that $A$
    is Gorenstein, but not strongly Gorenstein, as $3\ne \dim V$. 

Note that since $\kvg$ is integral over $A$, the height of $\Ann_A(g)$ is also
3, and as $A$ is Cohen-Macaualay, the depth of $\Ann_A(g)$ must be 3. Theorem
\ref{MainTheorem} implies that there exists a $p$-power $q$ such that $g^q =
0$. \done \end{eg}

We end with an example which shows that even in the non-modular case, the
good behaviour of separating algebras is not guaranteed by that of the
invariant ring.

\begin{eg}
Let $G=C_4$ be the cyclic group of order 4, and let $\kk$ be a field of odd characteristic containing a primitive fourth root of unity~$\zeta$. Consider the $2$-dimensional representation $V$ of $C_4=\langle \sigma\rangle$ given by\! \renewcommand{\arraystretch}{0.7}
\renewcommand{\arraycolsep}{2pt} $\sigma\mapsto {\tiny\left(\begin{array}{cc}\zeta&0\\0&\zeta\end{array}\right)}$. If $\kv=\kk[x,y]$, then $\kv^{C_4}=\kk[x^4,x^3y,x^2y^2,xy^3,y^4]$. On points where $x^4$ is zero, the function $x^2y^2$ also takes the value zero; on any other point $v$, we have $x^2y^2(v)=\frac{(x^3y(v))^2}{x^4(v)}$. Thus, the value of $x^2y^2$ as a function is entirely determined by the value of $x^4$ and $x^3y$. Therefore $x^4,x^3y,xy^3,y^4$ form a geometric separating set and $A:=\kk[x^4,x^3y,xy^3,y^4]$ is a geometric separating algebra. Note that $A$ is not Cohen-Macaulay, since the hsop $x^4,y^4$ does not form a regular sequence (despite being coprime) \cite[Exercise~2.1.18]{BrunsHerzog}. 

On points where $x^4$ is zero, $xy^3$ is zero, and on any other point $v$, $xy^3(v)=\frac{(x^3y(v))^3}{(x^4(v))^2}$. Thus, the hypersurface $\kk[x^4,x^3y,y^4]$ is also a graded geometric separating algebra.  Therefore, we have a Cohen-Macaulay graded geometric separating algebra, inside a non Cohen-Macaulay graded geometric separating algebra, inside a Cohen-Macaulay invariant ring.\done
\end{eg}

%%%%%%%%%%%%%%%%%%%%%%%%%%%%%%%%%%%%%%%%%%%%%%%%%%%%%%%%%%%%%%%%%%%%%%%%%%%%%%%%%%%%%%%%%%%%%%%%%%%%%%%%%%%%%%%%%%%%%%%%%%%%%%%%%%%%%%%%%%%%%%%%%%%%%%%%%%%%%%%%%%%%%%%%%%%%%%%%%%%%%%%%%%%%%%%%%%%%%%%%%%%%%%%%%%%%%%%%%%%%%%%%%%%%%%%%%%%%%%%%%%%%%%%%%%%%%%%%%%%%%%%%%%%%%%%%%%%%%%%%%%%%%%%%%%%%%%%%%%%%%%%%%%%%%%%%%%%%%%%%%%%%%%%%%%%%%%%%%%%%%%%%%%%%%%%%%%%%%%%%%%
%%%%%%%%%%%%%%%%%%%%%%%%%%%%%%%%%%%%%%%%%%%%%%%%%%%%%%%%%%%%%%

\bibliographystyle{plain}
\bibliography{MyBib}
\end{document}